\documentclass[12pt,a4paper]{amsart}
\usepackage{amssymb,amsmath,amsthm}
\usepackage{graphicx}
\usepackage{hyperref}
 \usepackage[color=green!40]{todonotes}
\usepackage{mathtools}
\usepackage{enumerate}

\newcommand\cA{{\mathcal A}}
\newcommand\cB{{\mathcal B}}

\newcommand\cF{{\mathcal F}}
\newcommand\cS{{\mathcal S}}

\newcommand\cG{{\mathcal G}}
\newcommand\cH{{\mathcal H}}

\newcommand\cQ{{\mathcal Q}}

\makeatletter
\newtheorem*{rep@theorem}{\rep@title}
\newcommand{\newreptheorem}[2]{%
\newenvironment{rep#1}[1]{%
 \def\rep@title{#2 \ref{##1}}%
 \begin{rep@theorem}}%
 {\end{rep@theorem}}}
\makeatother

\theoremstyle{plain}
\newtheorem{theorem}{Theorem}[section]
\newreptheorem{theorem}{Theorem}
\newtheorem{lemma}[theorem]{Lemma}

\newtheorem{proposition}[theorem]{Proposition}

\theoremstyle{definition}

\newcommand\cref[1]{Corollary~\ref{cor:#1}}

\title[Identification of a monotone Boolean function]{Identification of a monotone Boolean function\\ with $k$ ``reasons" as a combinatorial search problem}

\date{}
\author[]{D\'aniel Gerbner$^{a}$, András Imolay$^{b}$, Gyula O.H. Katona$^{a}$, D\'aniel T. Nagy$^{b}$, Kartal Nagy$^{b}$, Bal\'azs Patk\'os$^{a}$, Domonkos Stadler$^{b}$, Krist\'of Z\'olomy$^{b}$,  
\\
\small $^a$ HUN-REN Alfr\'ed R\'enyi Institute of Mathematics \\
\small $^b$ ELTE E\"otv\"os Lor\'and University} 
\begin{document}

\begin{abstract}
    We study the number of queries needed to identify a monotone Boolean function $f:\{0,1\}^n \rightarrow \{0,1\}$. A query consists of a 0-1-sequence, and the answer is the value of $f$ on that sequence. It is well-known that the number of queries needed is $\binom{n}{\lfloor n/2\rfloor}+\binom{n}{\lfloor n/2\rfloor+1}$ in general. Here we study a variant where $f$ has $k$ ``reasons'' to be 1, i.e., its disjunctive normal form has $k$ conjunctions if the redundant conjunctions are deleted. This problem is equivalent to identifying an upfamily in $2^{[n]}$ that has exactly $k$ minimal members. We find the asymptotics on the number of queries needed for fixed $k$. We also study the non-adaptive version of the problem, where the queries are asked at the same time, and  determine the exact number of queries for most values of $k$ and $n$.
\end{abstract}

\maketitle

\section{Introduction}

Let $[n]=\{ 1,2,\ldots, n\}$ be an $n$-element set (an \textit{$n$-set} for short). If $A\subseteq [n]$ is a subset, then its {\it characteristic vector} is $(x_1, x_2, \ldots , x_n)$ where $x_i=1 \; (1\leq i\leq n) $ if and only if $i\in A$, otherwise $x_i=0$.
A {\it Boolean function of} $n$ {\it variables} is a function $f(x_1, x_2, \ldots , x_n)$ giving the value $0$ or $1$ for every 0-1-sequence $(x_1, x_2, \ldots , x_n)$.
A Boolean function can be completely described by the family $\cH (f)$ of subsets of $[n]$ whose characteristic vectors give the value 1. General introductory texts on Boolean functions are \cite{CH,OD}.

Another way of describing Boolean functions uses the logical operations $\land, \lor$. 
Here $\land$ denotes ``and", 1 and 0 represent YES and NO, respectively.  That is, $x\land y=1$ if and only if both $x$ and $y$ are 1. Similarly, $\lor$ denotes ``or", that is $x\lor y=1$ if and only if at least one of $x$ and $y$ is 1. The {\it negation} $\bar{x}=1-x$ is also used. Any expression using these operations for the variables 
$x_1, x_2, \ldots , x_n$ defines a Boolean function. A {\it conjunction} is an expression of the following form:
$$x_{i_1}^{\prime}\land x_{i_2}^{\prime}\land \ldots \land x_{i_k}^{\prime}$$
where $x_{i_j}^{\prime}$ is either $x_{i_j}$ or $\overline{x_{i_j}}$ and $1\leq i_1 < i_2 < \ldots < i_k\leq n$. A {\it disjunctive normal form } is a disjunction  of some conjunctions.
In other words, several conjunctions are connected by the operation $\lor$.

The {\it identification} of a Boolean function (of $n$ variables) is an algorithm identifying the Boolean function, that is determining its value for all binary sequences of length $n$ knowing the value for some sequences. Of course, if no preliminary information is given on the function, then we cannot save anything, the value of the function should be queried for every sequence. A very natural assumption is  monotonicity. We say that the Boolean function 
$f(x_1, x_2, \ldots , x_n)$ is {\it monotone} if $f(x_1, x_2, \ldots , x_n)=1$ and 
$x_1\leq y_1, x_2\leq y_2, \ldots, x_n\leq y_n$ imply $f(y_1, y_2, \ldots , y_n)=1$.
In this case the family $\cH (f)$  of subsets giving value 1 has the following property: if $A\in \cH (f)$ and $A\subseteq B$ then $B\in \cH (f)$. In other words, for a monotone Boolean function the family $\cH (f)$ contains all the supersets of its members. Let us call such a family an {\it upfamily}. Finally, in the logical representation of the Boolean function, one can see that it is monotone if no negation appears in the disjunctive normal form, that is, it contains no $\overline{x_i}$. Conversely, every monotone Boolean function can be represented with a disjunctive normal form without negations.

There are two versions of identification  algorithms. A {\it non-adaptive algorithm} asks all the queries (questions, tests) at once and the unknown Boolean function should be determined based on the answers to these queries. The mathematical problem is to minimize the number of queries. In an {\it adaptive algorithm}, 
the choice of the next query may depend on  previous answers. In this case, we want to minimize the number of queries in the worst case.

It is easy to see that the condition that the function is monotone does not help in the case of non-adaptive algorithms.  We still have to ask all the values  to be able to identify the function. But the situation for adaptive algorithms is different. Korobkov \cite{Koro} proved that the adaptive algorithm needs at least
$$\binom{n}{\lfloor \frac{n}{2}\rfloor }+\binom{n}{\lfloor \frac{n}{2}\rfloor +1 }$$
queries in the worst case. On the other hand Hansel \cite{H} gave an algorithm with at most this number of queries. Selezneva and Liu \cite{SeleznevaLiu} showed that the same number of queries are needed if one of the aswers may be erroneous. 
Aslanyan and Sahakyan \cite{AS} generalized the methods for a more general function.
In \cite{AKS}, a special class of monotone Boolean functions was considered. They found a faster algorithm for this case.
Further algorithmic aspects of these problems are discussed in \cite{BorosHammer,Gainanov,Sok,Tono}. The average number of queries needed is studied in \cite{TT}.  

 In the present paper we assume that the monotone Boolean function to be determined has ``$k$ reasons" to be 1. More precisely, it can be represented as a disjunctive normal form consisting of $k$ conjunctions without negations. In terms of the family $\cH(f)$ this is equivalent to the condition that the number of minimal members is exactly $k$. Throughout the paper this representation will be used.

In the fundamental question of the theory of Combinatorial Search the $n$-set $[n]$ contains some ``wrong, defective" elements and they should be found by testing subsets  $A\subseteq [n]$. Let $D$ denote the set of defective elements. The result of the test is YES if $A\cap D\not= \emptyset$ and NO otherwise. There are infinitely many variants and generalizations, and in all such variants, the result of the test depends on some property of the intersection $A\cap D$. See the following books: \cite{AW, Aig, DH}. 

We know that $\cH(f)$ is an upfamily if $f$ is monotone. It is uniquely determined by its minimal members, so actually we need to find the family $\cS(f)$ of minimal members of $\cH(f)$. It is obvious that if $A, B\in \cS, A\not= B$ then $A\not\subseteq B$. Families having this property are called {\it antichains}.  The main differences between this problem and traditional Combinatorial Search are as follows: instead of a single subset $D$, we are looking for a $k$-member antichain. Queries involve choosing a subset $A \subseteq [n]$ to determine whether $A$ is a superset of at least one member of the unknown antichain. Thus, the problem reduces to identifying the hidden $k$-member antichain using as few such queries as possible.

Summarizing, our original problem of identification of a monotone Boolean function with $k$ reasons to be 1 is equivalent to the following problem of Combinatorial Search. An unknown $k$-member antichain is given in $[n]$. We can ask subsets of $[n]$ and the result of the query is that it is a superset of at least one member of the antichain, or it is not. We have to find a quick
algorithm using such queries uniquely determining the hidden antichain.

\subsection{Results}

As mentioned above, we will consider two versions of the problem of determining an antichain of $k$ sets. 
Let $f(n,k)$ denote the smallest number of queries that are enough to identify any antichain with $k$ members using an \textit{adaptive} algorithm. 

In the \textit{non-adaptive} version all the queries have to be asked in advance. Therefore, a non-adaptive algorithm corresponds to a family $\cG \subseteq 2^{[n]}$. We want to find a family $\cG$ that can identify  any antichain with $k$ members, i.e., for any two antichains $\cH$ and $\cH'$ with $|\cH|=|\cH'|=k$ there is a set $G\in\cG$ that is answered differently depending on whether $\cH$ or $\cH'$ is the antichain we want to find. It means that there is an $H\in\cH$ that is a subset of $G$ while no member of $\cH'$ is a subset of $G$, or the other way around. Let $h(n,k)$ denote the smallest size of a family $\cG$ with this property.

\begin{proposition}\label{egy}
    For any $n$ we have $f(n,1)=h(n,1)=n$.
\end{proposition}

\begin{proof}
    
The upper bound is given by the family consisting of the $(n-1)$-element sets. Observe that $x$ is contained in the unique element $H$ of $\cH$ if and only if $[n]\setminus \{x\}$ is answered NO, thus we can find each element of $H$. The lower bound is given by the so-called \textit{information theoretical lower bound}: there are $2^n$ possible antichains of size 1, and $q$ queries have $2^q$ possible sequences of answers. We have to be able to identify each possible antichain, thus distinct antichains have to result in distinct sequences of answers, showing that $2^n\le 2^q$.
\end{proof}

One of our main results is the almost complete determination of the function $h(n,k)$.

\begin{theorem}\label{nona} Let $n\ge 3$. Then
    \begin{displaymath}
h(n,k) = 
\left\{ \begin{array}{l l}
n & \text{if } k = 1, \vspace{1mm} \\
\sum_{i=n-k}^{n-1} \binom{n}{i} - 1 & \text{if } 1 < k < n, \vspace{1mm} \\
2^n-4 & \text{if } n \leq k \leq \binom{n-2}{\lfloor n/2\rfloor-1}+1.
\end{array}
\right.
\end{displaymath}
\end{theorem}

 We remark that $h(1,1)=1$ and $h(2,1)=2$ by Proposition~\ref{egy}, and $h(2,2)$=0, as there is only one antichain on $[2]$ with two subsets.

\medskip

Let us continue with the adaptive case. It turns out that, as we shall see, $f(n,k)$ is closely related to the number of minimal covers. Given a family $\cF$ of subsets of $[n]$, we say that a set $A\subseteq [n]$ is a \textit{cover} (or \textit{transversal}) of $\cF$ if $A\cap F\neq\emptyset$ for each $F\in \cF$. We say that $A$ is a \textit{minimal cover} of $\cF$ if $A$ is a cover of $\cF$ but proper subsets of $A$ are not covers of $\cF$. We let $MC(\cF)$ denote the family of minimal covers of  $\cF$. Let 

\[g(n,m):=\max\{|MC(\cF)|:\, \cF\subseteq 2^{[n]},\, |\cF|=m\}.\]
The next theorem establishes a connection between $f(n,k)$ and minimal covers.

\begin{theorem}\label{adap}
    \[ f(n,k)\le\sum_{m=1}^{k-1}g(n,m)+kn. \]
\end{theorem}

While the largest cardinality of the family of minimal covers is well-studied, see e.g. \cite{hy,lt}, we were unable to find any results on $g(n,m)$. We can determine $g(n,m)$ in the case $m$ is a constant and $n$ goes to infinity, and this, together with Theorem \ref{adap}, implies an upper bound on $f(n,m)$ which will turn out to be asymptotically sharp.

\begin{proposition}\label{fixadap}
    For any $m$ there exists $n_0=n_0(m)$ such that if $n\ge n_0$, then we have $g(n,m)=\prod_{i=0}^{m-1}\lfloor\frac{n+i}{m}\rfloor=(1+o(1))(\frac{n}{m})^m$ and the only families $\cF\subseteq 2^{[n]}$ of size $m$ with $|MC(\cF)|=g(n,m)$ are partitions with the property that $||F|-|F'||\le 1$ holds for all $F,F'\in \cF$. 
\end{proposition}

\begin{theorem}\label{adaptkbig}
    For any fixed $k\ge 3$ we have $f(n,k)=(1+o(1))\left(\frac{n}{k-1}\right)^{k-1}$.
\end{theorem} 
    
Finally, we deal with the missing case $k=2$.

\begin{proposition}\label{kketto}
    We have $f(n,2)=2n$ if $n\ge 12$.
\end{proposition}

Notation:  When the base set $[n]$ is fixed, we use $\overline{A}=[n]\setminus A$ to denote the complement of a set $A \subseteq [n]$. Furthermore, for a family $\cF$ we write $\overline{\cF}$ to denote $\{\overline{F}: F\in \cF\}$, the family of complements of sets in $\cF$.

\medskip

\textbf{Organization and terminology.} In Section 2, we shall prove our non-adaptive result, Theorem \ref{nona}, while Section 3 contains the proofs of all our adaptive results. 

We will use standard Combinatorial Search Theory terminology: the search problem will be considered as a two-player game between \textit{Questioner} and the \textit{Adversary}. An upper bound proof requires a strategy for Questioner which is just a set of queries in the non-adaptive case, and a binary tree with nodes labeled with queries and one edge leaving towards its child is labeled YES, while the other such edge is labeled NO.  In the non-adaptive case, a strategy solves the problem, if for any possible set of answers to the queries of Questioner, there is at most one antichain $\cS$ of $k$ sets that is consistent with the answers, i.e. for each query $Q$ of the Questioner there exists $H\in\cS$ with $H \subseteq Q$ if and only if the answer is YES. In the adaptive case, the same should hold for the set of query-labels and answer-labels along any branch of the binary tree.

For a lower bound proof we need a strategy for the Adversary. For any non-adaptive query set of size at most $M$, Adversary should be able to give answers with which at least two antichains of size $k$ are consistent. While in the adaptive case, the Adversary should have a strategy of answering queries such that at least one branch of the above binary tree is long if at the end there is only one consistent antichain. While answering, the Adversary might give away information for free (for example name the first $k-1$ sets of the antichain) that helps the Questioner to reduce the number of queries, but also helps us to analyze the number of queries needed.

\section{The non-adaptive case}

We prove Theorem \ref{nona} in this section, except for the case $k=1$, which is in Proposition \ref{egy}.

\medskip

\begin{proof}[Proof of Theorem~\ref{nona}] The case $k=1$ was addressed in Proposition \ref{egy}. 

\medskip

\textbf{Upper bound for the case $1<k<n$}

\smallskip

Consider the query family $\cQ$ consisting of all sets of size at least $n-k$ and at most $n-1$, except one of the $(n-1)$-element sets, say $[n-1]$. We have $|\cQ|=\sum_{i=n-k}^{n-1}\binom{n}{i}-1$. Assume that we have asked all sets in $\cQ$ and let $\cA$ be the family of sets $A$ with the property that each set $Q\supseteq A$ with $Q\in \cQ$ was answered YES. Let $\cA_0$ be the subfamily of minimal members of $\cA$. We will show that for the desired antichain $\cS=\{H_1,\ldots,H_k\}$, we have $\cS\subseteq \cA_0\subseteq \cS\cup \{[n-1]\}$. This would imply that we can determine $\cS$ as if $\cA_0$ has size $k$, then $\cS=\cA_0$, otherwise $\cA_0$ has size $k+1$, contains $[n-1]$, and then $\cS=\cA_0\setminus \{[n-1]\}$.

Clearly, each member of $\cS$ is in $\cA$ as for any query $Q\supseteq H_i$ the answer must be YES because of $H_i$. We next show $\cS\subseteq \cA_0$. Let $G\subsetneq H_i\in \cS$. We need to show that $G\notin \cA$, that is, a NO-answered query containing $G$ was asked. If $n-k\le |G|\le n-2$, then $G$ itself was asked and since $G\notin \cS$ ($\cS$ is an antichain!), we have the NO-answered query. As $k>1$, all members of $\cS$ have size at most $n-1$, so $G$ has size at most $n-2$. Finally, if $|G|<n-k$, then it is enough to show that there exists an $(n-k)$-set $K$ with $G\subseteq K$ and $H_i\not\subseteq K$ for all $i$, as then $K$ is a query with answer NO containing $G$. As $G$ is contained in some $H_i$, we have that for all $1 \le j \le k$ there exists $x_j\in H_j \setminus G$, as otherwise $H_j\subseteq G\subseteq H_i$, contradicting that $\cS$ is an antichain. 
Any $(n-k)$-subset of $[n]\setminus \{x_1,x_2,\ldots,x_k\}$ is a good choice for $K$ (note that the $x_j$'s may coincide).  This shows that $\cS\subseteq \cA_0$, indeed.

Assume that there exists $A\in \cA_0$ with $A\not\in\cS$. Then, as $\cA_0$ is an antichain and $\cS\subseteq \cA_0$, $A$ does not contain any $H_i$, thus there exists $x_i\in H_i \setminus A$ for each $1 \le i\le k$. Then $Q:=[n]\setminus\{x_1,\ldots,x_k\}$ contains $A$. If at least two of the $x_i$'s are distinct or all of them are the same, but not equal to $n$, then $Q$ is a query and thus, by $A\in \cA_0\subseteq \cA$, was answered YES, but does not contain any $H_i$, a contradiction.
If $x_1=x_2=\ldots=x_k=n$, then $n\in H_i$ for all $i$ and $n\notin A$. If $|A|\le n-2$, then there exists a query $Q$ with $A\subseteq Q\subsetneq [n-1]$. This is answered NO as all the sets $H_i$ contain $n$, so $A \notin \cA$. Therefore, the only possible extra member of $\cA_0$ compared to $\cS$ is $[n-1]$, as desired. Thus $h(n,k)\le |\cQ|=\sum_{i=n-k}^{n-1}\binom{n}{i}-1$.
This completes the proof of the upper bound.

\medskip

\textbf{Lower bound for the case $1<k<n$.}

\smallskip

Let us assume that we queried the family $\cA$ but did not query a set $A_0$ with $n-k\le |A_0|\le n-2$. Let $|A_0|=n-\ell$, $\overline{A_0}=\{w,w_1,\ldots, w_{\ell-1}\}$ and choose $x\in A_0$. 
Let $A_i=A_0\setminus \{x\}\cup \{w_i\}$ for $1\le i\le \ell-1$. For $\ell-1<i< k$, let $A_i$ be an $(n-\ell)$-set not contained in 
$A_0\cup \{w\}$, such that the sets $A_\ell, A_{\ell+1}, \dots, A_{k-1}$ are different.
Let us show the existence of the sets $A_i$ with $\ell-1<i< k$. We need $\binom{n}{n-\ell}-(n-\ell+1)-(\ell-1)=\binom{n}{n-\ell}-n\ge k-\ell$ as the left hand side is the number of $(n-\ell)$-sets with the subsets of $A_0\cup \{w\}$ and the $A_i$ ($1\le i\le \ell-1$) removed.
This inequality clearly holds if $\ell=k$, and otherwise the left-hand side is at least $n-2$, while the right-hand side is at most $n-2$. 

Now the Adversary uses the following strategy: he answers YES to each set 
that contains any of the sets $A_i$ with $0\le i\le k-1$ 
and the answer is NO to all other queries. 
Then it is possible that $\cS=\{A_0,\ldots,A_{k-1}\}$, and it is also possible that $\cS=\{A_0\cup\{w\}, A_1,A_2,\ldots,A_{k-1}\}$. Therefore the family $\cA$ cannot determine $\cS$, so one indeed needs to query all sets of size between $n-k$ and $n-2$.

We claim further that there cannot be two unqueried sets $A_0, A_1$ with $|A_0| = |A_1| = n-1$.
Let $A_2, A_3, \ldots, A_{k}$ be distinct sets with $|A_i|=n-1$ different from $A_0$, $A_1$. Clearly, there are $n-2 \geq k-1$ such sets.
Assume that the Adversary answers YES to each query containing any of $A_0, A_1, \ldots, A_{k}$ and NO to all other queries. Then $\mathcal{S}= \{A_0, A_2, A_3, \ldots, A_{k}\}$ and $\mathcal{S}'=\{ A_1, A_2, \ldots, A_{k}\}$ are both still possible, which is a contradiction. So to determine $\cS$, we indeed need to query all but at most one set of size $n-1$.

In conclusion, every set of size between $n-k$ and $n-1$ must be queried, except for at most one, and so $h(n,k)\ge \sum_{i=n-k}^{n-1}\binom{n}{i}-1$ as claimed.

\medskip

\textbf{Upper bound for the case $k\ge n$.} 

\smallskip

We ask all sets except for $\emptyset$, $[n]$, $\{n\}$ and $[n-1]$, a total of $2^n-4$ sets. If the answers allow us to determine $\cS$, then $h(n,k)\le 2^n-4$ as claimed. 

Assume that we have asked these sets. If a set $A$ of size $m$ and all of its subsets of size $m-1$ were asked, then $A\in \cH$ if and only if each of those subsets is answered NO and $A$ is answered YES. This shows that we know each element of $\cS$ except possibly $\{n\}$, the 2-element sets containing $n$ and $[n-1]$. Note that $\emptyset$ and $[n]$ cannot be in $\cS$ since $k>1$. 

Observe that $[n-1]$ is not in $\cH$ if any of its $(n-2)$-element subsets are answered YES, thus if we have uncertainty here, then those $(n-2)$-element subsets are each answered NO and $[n-1]$ can be the only element of $\cS$ that does not contain $n$. In this case, $\{n\}$ cannot be in $\cH$ as $k \geq n \geq 3$. This implies that $\{n\}$ would have been answered NO, thus we know which 2-element sets are in $\cS$ as well. We know each element of $\cS$ except for $[n-1]$, thus we only have to check whether there are $k$ or $k-1$ other sets that we have already identified to be in $\cS$.

By the above paragraph, we can assume that some $(n-2)$-element subset of $[n-1]$ were answered YES, thus $[n-1]\not\in\cS$. If any $2$-element sets that contain $\{n\}$ are answered NO, then we can consider $\{n\}$ as being answered NO and we are done. Otherwise, either $\{\{1,n\},\{2,n\},\ldots,\{n-1,n\}\}$ or $\{n\}$ belong to $\cS$. We can simply count the other sets we already know to be in $\cS$ and check whether there are $k-1$ or $k-n+1$ of them.

\medskip

\textbf{Lower bound of the case $n\le k\le \binom{n-2}{\lfloor n/2\rfloor -1}+1$.} 

\smallskip

Let us assume that we queried the family $\cA$ but did not query a set $A_0$ with $2\le |A_0|=n-\ell\le n-2$. Without loss of generality, $[2]\subseteq A_0=[n-\ell]\subseteq [n-2]$. Let us pick a family $\cB$ of $k-\ell+1$ sets of size $\lfloor n/2\rfloor$ that each contains $n$ and does not contain $1$. As long as $k\le \binom{n-2}{\lfloor n/2\rfloor -1}+1$ 
we can pick such a family $\cB$.
Let $A_i=A_0\setminus \{2\}\cup \{n-i\}$ for $1\le i\le \ell-2$. Observe that $\cB\cup \{A_i:0\le i\le \ell-2\}$ is an antichain and so is $\cB\cup \{A_i: 1\le i\le \ell-2\}\cup \{A_0\cup \{n-\ell+1\}\}$.

Suppose the Adversary gives answers to every possible set in $2^{[n]}\setminus \{A_0\}$ as follows: YES to every set containing $A_0$ or $A_i$ or a set in $\cB$ and NO to other sets. Now we know that $A_i$ for $1\le i\le \ell-1$ is in $\cH$, and $\cB\subseteq \cS$. Additionally, either $A_0$ is in $\cS$, or $A_0\cup \{n-\ell+1\}$ is in $\cS$, showing that we cannot identify $\cS$.

Let us assume that we have not queried two sets of size 1, say $\{1\}$ and $\{2\}$. Then we pick a family $\cB$ of $k-n+1$ sets of size $\lfloor n/2\rfloor-1$ such that $B \cap [2]=\emptyset$ for all $B \in \cB$. We can pick such $\cB$ as long as $k\le \binom{n-2}{\lfloor n/2\rfloor-1}+n-1$. Suppose the Adversary answers YES to each set that contains 1 or 2 or a set in $\cB$, and NO to all other sets. Now it is possible that in addition to $\cB$, $\cS$ consists of $\{1\}$ and the sets of the form $\{2,i\}$, $3\le i\le n$, and it is also possible that in addition to $\cB$, $\cS$ consists of $\{2\}$ and the sets of the form $\{1,i\}$, $3\le i\le n$. (Both of these families are antichains.) This shows that we cannot identify $\cS$.

Finally, let us assume that we have not queried two sets of size $n-1$, say $[n-1]$ and $[n]\setminus \{n-1\}$. Let us pick a family $\cB$ of $k-1$ sets of size $\lfloor n/2\rfloor$+1 that each contains both $n-1$ and $n$.
We can pick such $\cB$ as long as $k\le \binom{n-2}{\lfloor n/2\rfloor-1}+1$. Suppose the Adversary answers YES to each query that contains at least one of the sets in $\cB$, and all other queries are answered NO. Then it is possible that $\cS$ consists of $[n-1]$ and $\cB$, and also possible that $\cS$ consists of $[n]\setminus \{n-1\}$ and $\cB$, showing that we cannot identify $\cS$. 

Summarizing the above: to determine $\cS$, one has to ask all subsets of $[n$] but at most one subset of size 1, at most one set of size $n-1$, and additionally one does not need to query $\emptyset$ and $[n]$. So all other $2^{n}-4$ subsets must be queries. This completes the proof of this case.
\end{proof}

\section{The adaptive case}

\begin{proof}[Proof of Theorem \ref{adap}]

We shall use Gainanov's algorithm \cite{Gainanov} that finds a member of $\cS$ in a set $A$ using $|A|$ queries, if we know that $A$ contains at least one member of $\cS$ (either because $A=[n]$ and $k>0$, or because $A$ was a query and was answered YES). Let $a=|A|$ and without loss of generality $A=[a]$. Let $\cS_0=\cH\cap 2^{A}$, then $\cS_0$ is non-empty.

     We start by asking $[a-1]$. If the answer is YES, then we know that there is a member of $\cS_0$ that does not contain the element $a$, while if the answer is NO, then we know that each member of $\cS_0$ contains the element $a$. In the first case the further queries will not contain the element $a$, in the second case they will. This way we continue with a non-empty subfamily $\cS_1$ of $\cS_0$: in the case of YES answer, $\cS_1$ consists of the sets in $\cS_0$ that do not contain the element $a$, in the case of NO answer, $\cS_1$ consists of the sets in $\cS_0$ that contain the element $a$ (which are all the sets in $\cS_0$).
    
    Then we remove the element $a-1$, i.e., we either ask $[a-2]$ or $[a-2]\cup\{a\}$. If the answer is YES, we know that the subfamily $\cS_2$ of $\cS_1$ consisting of sets that do not contain the element $a-1$ is non-empty, if the answer is NO, we know that  $\cS_2:=\cS_1$ consists of sets that contain the element $a-1$. We continue this way, and after $a$ steps we know each element of the only member of $\cS_a$.

   Using the above algorithm, first we find one member of $\cS$ using at most $n$ queries. 
Assume that we already identified a subfamily $\cS_m' \subseteq \cS$ with $|\cS_m'|=m$ for some $1 \leq m<k$.
   We show that it is possible to find one more member of $\cS$ with at most $g(n,m)+n$ queries.
We ask each set in $\overline{MC(\cS_m')}$. Observe that if $H\in \cS\setminus \cS_m'$, then $H$ does not contain any set in $\cS_m'$, hence for each $H'\in\cS_m'$ there is an element $v\in H'$ such that $v\not\in H$. Taking these elements for each set in $\cS'_m$, we obtain a set $H''$ that is a cover of $\cS_m'$ and is disjoint with $H$. Therefore, there is a minimal cover $H'''$ of $\cS_m'$ that is disjoint with $H$. The complement of $H'''$ contains $H$, thus is answered YES. We obtained a YES answer, thus we can find another member of $\cS$ using at most $n$ queries. Altogether we used at most $g(n,m)+n$ queries to find the $m+1$st element. Summing up we obtain the bound, completing the proof.
\end{proof}

\smallskip

\begin{proof}[Proof of Proposition \ref{fixadap}]
    The lower bound is given by the family of the $m$ parts of an equipartition of $[n]$, i.e. pairwise disjoint sets $H_1,H_2,\ldots,H_m$ with $\cup_{i=1}^mH_i=[n]$ and $||H_i|-|H_j||\le 1$ for all $1\le i< j\le m$.
    
    For the upper bound, let $|\cF|=m$. If $\cF$ consists of pairwise disjoint sets, then $|MC(\cF)|=\prod_{F\in \cF}|F|$ is maximized when $\cF$ is a partition and for all $F, F' \in \cF$, $||F|-|F'||\le 1$. If $\cF$ is not a family of pairwise disjoint sets, then let $A$ be the set of elements that are contained in more than 1 member of $\cF$ and let us write $a=|A|\ge 1$. Let $a_i$ be the number of elements of the $i$th set in $\cF$ in $[n]\setminus A$, so $a+\sum_{i=1}^ma_i\le n$. After renumbering, we may assume $a_1\ge a_2\ge \ldots \ge a_m$. The number of sets in $MC(\cF)$ avoiding $A$ is $\prod_{i=1}^m a_i$. As every  member of $MC(\cF)$ intersecting $A$ has size at most $m-1$, their number can be counted the following way. We pick a subset $A'$ of $A$ of size $j\le m-1$, then we pick one element from $F_i\setminus A$ for each $i$ with $A'\cap F_i=\emptyset$. Note that for any $A'$ there is at most $m-1-|A'|$ such $i$.  The resulting set is a cover but not necessarily minimal. However, each minimal cover intersecting $A$ is obtained this way. Hence
    \begin{equation}\label{eq}
    |MC(\cF)|\le \prod_{i=1}^ma_i+\sum_{j=1}^{\min\{a,m-1\}}\binom{a}{j}\prod_{i=1}^{m-j-1}a_i\le \prod_{i=0}^{m-1}\left\lfloor\frac{n-a+i}{m}\right\rfloor+man^{m-2},    
    \end{equation}
    where we used $a,a_i\le n$, $a\ge 1$ and thus $\binom{a}{j}\prod_{i=1}^{m-j-1}a_i\le a^jn^{m-j-1}\le an^{m-2}$ for any $1\le j\le \min\{a,m-1\}$.

    In view of (\ref{eq}), to prove the proposition, we need to bound $\prod_{i=0}^{m-1}\lfloor\frac{n+i}{m}\rfloor-\prod_{i=0}^{m-1}\lfloor\frac{n-a+i}{m}\rfloor$ from below, more precisely we need $\prod_{i=0}^{m-1}\lfloor\frac{n+i}{m}\rfloor-\prod_{i=0}^{m-1}\lfloor\frac{n-a+i}{m}\rfloor>man^{m-2}$. Observe first that this difference is monotone increasing in $a$. Therefore,  if $a\ge n/2$, then $\prod_{i=0}^{m-1}\lfloor\frac{n+i}{m}\rfloor-\prod_{i=0}^{m-1}\lfloor\frac{n-a+i}{m}\rfloor\ge \prod_{i=0}^{m-1}\lfloor\frac{n+i}{m}\rfloor-\prod_{i=0}^{m-1}\lfloor\frac{n/2+i}{m}\rfloor =(1-\frac{1}{2^m}+o(1)) (\frac{n}{m})^{m}$ which is larger than $man^{m-2}$ if $n$ is large enough. 

    Finally, note that $\prod_{i=0}^{m-1}\lfloor \frac{b+i}{m} \rfloor$ is the number of minimum covers of an equipartition of a $b$-set $B$ into $m$ parts. 
    Take two sets $B' \subsetneq B$ with cardinalities $b'$ and $b$ respectively, and fix an equipartition of $B$ such that its restriction to $B'$ is an equipartition of $B'$. 
    The expression $\prod_{i=0}^{m-1}\lfloor \frac{b+i}{m} \rfloor - \prod_{i=0}^{m-1}\lfloor \frac{b'+i}{m} \rfloor$
    counts those minimal covers of the equipartition of $B$ which are not contained in $B'$, including 
    those that contain one element $x$ from $B\setminus B'$ and one element from each part of $B'$ of which the extension in $B$ does not contain $x$. If $a\le n/2$, then the number of elements in the difference is $a$ and in each part of $B$, there remain at least half  of the original elements in the part of $ B'$, thus at least $\lfloor\frac{n}{2m}\rfloor$ elements. We obtain $\prod_{i=0}^{m-1}\lfloor\frac{n+i}{m}\rfloor-\prod_{i=0}^{m-1}\lfloor\frac{n-a+i}{m}\rfloor\ge a\lfloor \frac{n}{2m}\rfloor^{m-1}$ which is larger than $man^{m-2}$ if $n$ is large enough.
    This completes the proof. 
\end{proof}

\smallskip

\begin{proof}[Proof of Theorem \ref{adaptkbig}]
    Theorem \ref{adap} and Proposition \ref{fixadap} imply $f(n,k)\le \sum_{m=1}^{k-1}g(n,m)+kn\le kn+\sum_{m=1}^{k-1}(1+o(1))(\frac{n}{m})^m\le (1+o(1))(\frac{n}{k-1})^{k-1}$. This proves the upper bound.

    To see the lower bound, we describe the following simple strategy of the Adversary for $n>k$. He tells the Questioner that $\cS$ contains a partition $H_1,H_2,\ldots,H_{k-1}$ of $[n]$ with $|H_i|=\lfloor \frac{n+i-1}{k-1}\rfloor$ and the last member of $\cS$ is a set $H_k\in \overline{MC(\cS')}$ with $\cS'=\{H_1,H_2,\ldots,H_{k-1}\}$. So Questioner's only task is to determine $H_k$. 

       If we query a set that contains an $H_i$, then the answer is clearly YES, thus there is no point of such a query, so we assume that Questioner's query sets do not contain any member of $\cS'$.
  If we query a set that contains neither a member of $\overline{MC(\cS')}$ nor a member of $\cS'$, then the answer is clearly NO, 
    thus there is no point of such a query. If we query a set $Q$ that is a proper superset of a member of $\overline{MC(\cS')}$, then $\overline{Q}$ is a proper subset of a member of $MC(\cS')$, hence there is a set $H_i \in \cS'$ with $\overline{Q} \cap H_i=\emptyset$. Consequently, $H_i\subsetneq Q$ for some $i$. Therefore, the answer to $Q$ is YES and there is no point of such a query.

 We obtained that each query has to be a member of $\overline{MC(\cS')}$. As $n>k$, it is easy to see that $\cS' \cup \overline{MC(\cS')}$ is an antichain. Hence any member of $\overline{MC(\cS')}$ can be the unique set of $\cS \setminus \cS'$. The Adversary can answer NO to the first $\prod_{i=1}^{k-1}\lfloor \frac{n+i-1}{k-1}\rfloor-2$ such queries. Then all we know is that among the remaining $2$ members of $\overline{MC(\cS')}$, one of them belongs to $\cS$ and the other does not, but we have no way to find out which one belongs to $\cS$ (except for asking one more query of course). This proves the lower bound.
\end{proof}

Before proving Proposition \ref{kketto}, we determine the number of antichains of size 2. This is a well-known exercise, we include its proof for completeness.

\begin{lemma}
    The number of antichains of size 2 in $[n]$ is $A(n):=2^{2n-1}-3^n+2^{n-1}$.
\end{lemma}

\begin{proof}
    First we show that the number of pairs $\{A,B\}$ with $A\subseteq B\subseteq [n]$ is $3^n$. We identify such pairs with vectors of length $n$ and values $0,1,2$. If coordinate $i$ has value 0, then $i$ belongs to none of $A$ and $B$, if coordinate $i$ has value 1, than $i$ belongs to $B$ but not $A$, and if coordinate $i$ has value 2, then $i$ belongs to both $A$ and $B$. Clearly we obtain such a pair for every vector, and for every pair we can define such a vector by reversing this. This implies that the number of chains of size 2 is $3^n-2^n$.

    There are $2^{n-1}(2^n-1)=2^{2n-1}-2^{n-1}$ pairs of distinct sets in $[n]$, $3^n-2^n$ of those are chains, the remaining $2^{2n-1}-2^{n-1}-3^n+2^n$ pairs are antichains, completing the proof.
\end{proof}

We remark that using this lemma, the information theoretical lower bound $f(n,2)\ge \log_2 A(n)$ yields $f(n,2)\ge 2n-1$ for $n \geq 5$. In the  proof of the lower bound of Proposition \ref{kketto}, we will use the following simple idea. If an algorithm used $m$ queries and there are at least $M$ antichains of size 2 that are consistent with the answers to the $m$ queries asked, then
the $M$ antichains of size 2 must have distinct sequences of answers. Therefore,
at least $\lceil\log_2M\rceil$ further queries are needed, hence the total number of queries is at least $m+\lceil\log_2M\rceil$. We will only have to consider $m=1$.

\smallskip

\begin{proof}[Proof of Proposition \ref{kketto}] The upper bound uses Gainanov's algorithm, see the proof of Theorem \ref{adap}. First we find a member $H$ of $\cS$ using $n$ queries. Then we can find the second member also, using $n$ queries the following way. We first ask $[n]\setminus \{h_1\}$ for some $h_1\in H$. If the answer is YES, we continue inside $[n]\setminus \{h_1\}$ and finish using $n-1$ further queries. If the answer is NO, we know that the other member $H'$ of $\cS$ contains $h_1$ and ask $[n]\setminus \{h_2\}$ for some $h_2\in H$. We continue this way, always picking a new vertex of $H$. For some $i\le |H|$, we obtain a YES answer, otherwise $H'$ contains $H$. Then we know that $h_1,\ldots,h_{i-1}\in H'$, thus they are contained in each further queries, while $h_i\not \in H'$. We continue by removing a potential vertex of $H'$ from the query set as in the proof of Theorem \ref{adap}. Each answer tells us whether that element is in $H'$ or not, thus with $n-i$ further queries we identify $H'$.

For the lower bound, consider the first query $A$. Assume first that $|A|=n-1$. Then an answer NO means that both members $H_1$ and $H_2$ of $\cS$ contain the element $x=[n]\setminus A$, thus $H_1\setminus \{x\}$ and $H_2\setminus \{x\}$ is an antichain. This implies that there are $A(n-1)=2^{2n-3}-3^{n-1}+2^{n-2}$ antichains where the answer is NO. Therefore, there are \[A(n)-A(n-1)=2^{2n-1}-3^n+2^{n-1}-2^{2n-3}+3^{n-1}-2^{n-2}>2^{2n-2}\] antichains if the answer is YES. Here we use that $n\ge 6$. By the argument before the proof, more than $2n-2$ further queries are needed to find the antichain if the answer was YES to a query of size $n-1$. So a total of at least $1+2n-1=2n$ queries are required.

Assume now that $|A|\le n-2$. Then an answer YES means that one of the members of $\cS$ is in $A$. There are $A(n-2)=2^{2n-5}-3^{n-2}+2^{n-3}$ antichains with $2$ members inside $A$. For the number of antichains with one member inside $A$ and the other member not inside $A$, we use the simple upper bound $2^{n-2}(2^n-2^{n-2})$. This counts each pair $(A',B)$ with $A'\subseteq A$, $B\not\subseteq A$, including those pairs where $A'$ is a subset of $B$. This is enough for us, since it means that in the case of a NO answer to $A$, there are at least \[2^{2n-1}-3^{n}+2^{n-1}-(2^{2n-5}-3^{n-2}+2^{n-3})-2^{n-2}(2^n-2^{n-2})>2^{2n-2}\] antichains. Here we use that $n\ge 12$. Again, the argument from the paragraph before the proof completes the proof.
\end{proof}
\smallskip
\textbf{Acknowledgement}: Research was carried out at the Combinatorial Search Theory Seminar of the Alfr\'ed R\'enyi Institute of Mathematics.

\smallskip
\textbf{Funding}: Research was supported by the National Research, Development and Innovation Office - NKFIH under the grants FK 132060, KKP-133819, SSN135643 and K132696. András Imolay was supported by the EKÖP-24 University Excellence Scholarship Program of the Ministry for Culture and Innovation from the source of the National Research, Development and Innovation Fund. Kristóf Zólomy was supported by the Thematic Excellence Program TKP2021-NKTA-62 of the National Research, Development and Innovation Office.

\end{document}